\documentclass{article}
\usepackage{amssymb}

\usepackage{amsmath}
\usepackage{chicago}

\newtheorem{theorem}{Theorem}[section]

\newtheorem{corollary}{Corollary}

\newtheorem{proposition}{Proposition}[section]

\input{tcilatex}

\begin{document}

\title{Optimal Convergence Trading}
\author{Vladislav Kargin\thanks{%
Cornerstone Research, 599 Lexington Avenue, New York, NY 10022; slava@bu.edu}%
}
\maketitle

\begin{abstract}
This article examines arbitrage investment in a mispriced asset when the
mispricing follows the Ornstein-Uhlenbeck process and a credit-constrained
investor maximizes a generalization of the Kelly criterion. The optimal
differentiable and threshold policies are derived. The optimal
differentiable policy is linear with respect to mispricing and risk-free in
the long run. The optimal threshold policy calls for investing immediately
when the mispricing is greater than zero with the investment amount
inversely proportional to the risk aversion parameter. The investment is
risky even in the long run. The results are consistent with the belief that
credit-constrained arbitrageurs should be risk-neutral if they are to engage
in convergence trading.
\end{abstract}

\begin{quotation}
\bigskip 

\textsl{Myron [Scholes] once told me they are sucking up nickels from all
over the world. But because they are so leveraged that amounts to a lot of
money.}

\textsl{Merton Miller about the essence of arbitrage.}
\end{quotation}

\section{Introduction}

Arbitrageurs are people who detect inconsistencies in asset prices and
invest in them hoping that the inconsistencies will be eliminated. The
waiting time is often uncertain and since the arbitrageur depends on the
willingness of other people to lend him money, the irrationality of
creditors may lead to great debacles long before prices converge to
consistent values. The notorious story of the arbitrage fund LTCM that lost
90 percent of its value on ``riskless''\ deals illustrates the importance of
credit constraints. So what policy should the arbitrageur pursue when
creditors impose borrowing constraints? In particular, can the arbitrageur
allocate the available funds in such a way as to eliminate all the long-run
risk?

If this risk elimination is possible, the mispricings should be equally
attractive to risk-averse as well as risk-neutral investors. However, a
popular view asserts that arbitrageurs, unlike other investors, should be
risk-neutral. Is there any ground for this belief? The present article
provides a justification by solving for arbitrageurs' optimal policies under
several types of constraints and showing that under some of them the
long-run risk cannot be eliminated.\textbf{\ }Risk-averse investors that
face those constraints are not interested in small mispricings.

The paper investigates two classes of constraints that lead to strikingly
different results. Under constraints from the first class, the arbitrageur
can only change the leverage slowly. In practice, borrowing additional funds
takes time: The arbitrageur must apply for new credit, provide an
explanation why he needs it and wait for a decision. Depending on the
situation, the process might take from several minutes to several days. In
addition a rapid increase in a position adversely affects prices, so in
their own interests arbitrageurs must accumulate positions slowly.

For this class of policies, the main result is that the optimal policy is
linear in the mispricing, and independent of the coefficient of risk
aversion. The variance of the portfolio wealth does not grow with time.
Thus, under this constraint the long-run risk can be expunged.

Constraints of the second type are stronger: The arbitrageur cannot change
leverage except by closing the position. The motivation is that the
arbitrageur is often restricted in his ability to change the leverage --
even if the need arises. Higher leverage is mostly needed when mispricing is
increasing and the investment account shows negative performance.
Unfortunately, this is the worst time to ask for new credit because the
creditors hate to invest in accounts with negative performance. As Mark
Twain said: ``A banker is a fellow who lends you his umbrella when the sun
is shining but wants it back the minute it rains.''

For policies in this class, the main result is that the long-run risk cannot
be completely removed. Consequently, the arbitrageur will invest an amount
that is inversely proportional to his risk aversion.

These two examples suggest that what makes the convergence arbitrage risky
in the long run is the inability of the arbitrageur to change the investment
amount after the investment is committed. In particular, the results of the
second example are consistent with the belief that arbitrageurs should be
risk-neutral if they are to engage in convergence trading.

The results of the present paper match closely with results of %
\shortciteN{grossman_vila92}, who study the dynamic investment under a
constraint on investment amount. They find that the constraint essentially
makes the investor behave as if he were more risk-averse than he actually
is. Unlike in the present paper, however, the asset process is not
mean-reverting in \shortciteN{grossman_vila92}, so the investor could not
hope to eliminate the risk completely. Also the constraint is not exogenous
as in the present paper but a function of the investor's wealth. Because of
these differences it is difficult to conclude whether the similarity of
results is incidental or not. Both papers, however, support the view that
certain constraints increase long-run riskiness of investment projects.

In a recent paper about convergence trading, \shortciteN{longstaff_liu00}
use the Brownian bridge to model the mispricing process, an assumption on
the process that requires a fixed horizon at which the mispricing will be
effaced. By the nature of their model, they cannot draw conclusions about
long-run risks but they do find that arbitrageurs sometimes cannot eliminate
all risk at the end of the arbitrage period. This result is consistent with
results of the present paper.

Closely related is the literature on optimal dynamic investments with risky
assets. The seminal contributions to this literature were made in %
\shortciteN{samuelson69}, \shortciteN{merton69}, \shortciteN{merton71} and %
\shortciteN{merton73a}. Optimal investment in assets that follow a diffusion
process with mean-reverting returns was analyzed in \shortciteN{kim_omberg96}%
, \shortciteN{brennan_schwartz_lagnado97}, \shortciteN{campbell_viceira99}, %
\shortciteN{barberis00}, \shortciteN{wachter02}, %
\shortciteN{campbell_chan_viceira03}. The focus of the present paper is not
on general risky investments but on the optimal extraction of profit from
near arbitrage opportunities. Concequently, the paper comes to more definite
conclusions by using a generalization of the Kelly investment criterion,
which emphasizes the long-run behavior of portfolios and especially suitable
to modelling objectives of large institutional traders.

In addition, this paper computes the optimal leverage using a new method.
While Kim and Omberg ingeniously solve dynamic programming equations by
reducing them to a system of non-linear ordinary differential equations,
Campbell and Viceira derive an approximate solution by linearization of
Euler equations, and Wachter uses the martingale method of %
\shortciteN{cox_huang89} to separate consumption and financing decisions,
the present paper derives the solution by methods of stochastic control,
taking the advantage of the asymptotic investment criteria.

The rest of the paper is organized as follows. Section 2 describes the
model. Sections 3 and 4 derive the optimal differentiable and threshold
policies and describe their properties. Section 5 compares results for
differentiable and threshold policies and concludes.

\section{Model}

An investor can invest in a mispriced asset whose mispricing is measured by $%
x=\ln (p_{1}/p_{2})$. Here $p_{1}$ and $p_{2}$ are the asset's actual and
``correct'' prices, respectively. Mispricing follows the Ornstein-Uhlenbeck
process:%
\begin{equation}
dx_{t}=-\alpha x_{t}dt+\sigma dz_{t},
\end{equation}%
where $x_{t}$ is mispricing at time $t$, $z_{t}$ is a standard Wiener
process, $\sigma >0$ and $\alpha >0$. Parameter $\alpha $ measures the speed
of reversion to the correct price: The higher $\alpha $ is, the faster
mispricing $x$ drifts towards zero. Parameter $\sigma $ measures the size of
new mispricing shocks introduced into the process. It is also useful to
define $\Sigma =\sigma ^{2}/(2\alpha ),$ which is the variance of $x_{t}$ in
the long run.

Changes in mispricing induce changes in the arbitrageur's wealth through his
choice of the leverage coefficient $f(x)$: The change in the logarithm of
wealth is the product of the leverage coefficient and the change in the
mispricing,

\begin{equation}
du=f(x)dx.
\end{equation}%
Intuitively, a $1\%$ change in the mispricing results in a $f(x)\%$ change
in the investor's wealth. Later we will impose certain restriction on the
available leverage.

The arbitrageur's utility is a linear combination of the growth rates in the
expectation and variance of the portfolio wealth:%
\begin{equation}
U=\underset{t\rightarrow \infty }{\lim \inf }\frac{1}{t}\left[
E(u_{t})-\gamma \mathrm{Var}(u_{t})\right] ,
\end{equation}%
where parameter $\gamma $ measures risk aversion of the investor. The
optimization problem is to choose the leverage function $f(x)$ so that
utility $U$ is maximized.

What is the meaning of this maximization criterion? If $\gamma $ is $0,$
then the criterion is the same as the criterion of maximizing the
portfolio's long-run growth rate -- the Kelly criterion. When $\gamma >0,$
it introduces an additional term penalizing deviations from the expected
growth rate. This additional term assures the investor that maximizing $U$
protects him against the large deviations in the realized growth of his
portfolio from the expectations.

An example may perhaps add some insight into the investment criterion.
Suppose that the wealth of the investor follows a geometric Brownian motion
with constant parameters $\mu $ and $\sigma .$ Then the utility of the
investor is 
\begin{equation}
U=\mu -\gamma \sigma ^{2}.
\end{equation}%
This expression shows that the utility depends only on the parameters of the
process and on risk aversion but not on the investment horizon.

Another way to get an insight into this criterion is to compare it with the
objective under the classical single period Markowitz model. In the
Markowitz model the investor maximizes a linear combination of the
expectation and the variance of the portfolio return. Therefore, the present
model generalizes the Markowitz model to the dynamic setting by substituting
the expectation and variance of the single period return with the asymptotic
rates of increase in expectation and variance of the investor's wealth.

An important assumption that we adopt in this generalization is that the
investor is concerned only with long-run consequences of his policy. This
assumption simplifies the analysis considerably and appears to be realistic
for small investments by large institutions. In using this assumption we
follow \shortciteN{bielecki_pliska99} and %
\shortciteN{bielecki_pliska_sherris00}, who applied it to the analysis of
continuous investment policies in a similar situation.

\section{Optimal Differentiable Policy}

\label{lever:OptimalDifferentiableStrategy}

This section is about differentiable policies $f(x)$, for which 
\begin{align}
& f\in C^{1}(-\infty ,+\infty ),\text{ and} \\
|& f^{\prime }(x)|\leqslant K.
\end{align}

\noindent The policies from this class will be called $\mathbf{D-}$policies$%
. $ This class excludes policies that prescribe extremely rapid growth of
leverage with respect to mispricing. The following theorem is a cardinal
ingredient in showing that optimal $\mathbf{D}-$policies are linear.

\begin{theorem}
\label{lever:proposition1} Linear investment policies are the only $\mathbf{D%
}-$policies such that the variance of the logarithm of wealth $u_{t}$ is
asymptotically constant.
\end{theorem}

\noindent The proof uses a convenient representation for $u$: Let 
\begin{equation}
g(\xi )=:\int_{0}^{\xi }f(\zeta )d\zeta .
\end{equation}%
Then it is easy to check that 
\begin{equation}
u_{t}=u_{0}+g(x_{t})-g(x_{0})-\frac{\sigma ^{2}}{2}\int_{0}^{t}f^{\prime
}(x_{\tau })d\tau .  \label{representation}
\end{equation}

The intuition behind this representation is simple. The investor can
increase his wealth only if he increases his leverage when the mispricing
increases. The derivative $f^{\prime }(x)$ measures the sensitivity of the
leverage policy to mispricing, and (\ref{representation}) shows that the
change in the logarithm of wealth equals a multiple of the integral of $%
f^{\prime }(x)$ plus a stationary process. The more sensitive the leverage
policy is to mispricing, the greater the increase in the wealth induced by
local variations in mispricing. The addition of $g(x_{t})-g(x_{0})$ reflects
dependence of the wealth on the initial and final conditions.

\noindent \textbf{Proof of Theorem \ref{lever:proposition1}:} Taking the
variance of $u_{t}-u_{0}$ in (\ref{representation}) gives

\begin{equation}
\mathrm{Var}(u_{t}-u_{0})=\mathrm{Var}\left( g(x_{t})\right) -\sigma ^{2}%
\mathrm{Cov}\left( g(x_{t}),\int_{0}^{t}f^{\prime }(x_{\tau })d\tau \right) +%
\frac{\sigma ^{4}}{4}\mathrm{Var}\left( \int_{0}^{t}f^{\prime }(x_{\tau
})d\tau \right) .
\end{equation}

\noindent As $t$ increases, all terms except possibly the third one tend to
a finite limit. So, asymptotically, 
\begin{equation}
\mathrm{Var}(u_{t})\sim \mathrm{const}+\frac{\sigma ^{4}}{4}rt,
\end{equation}%
where 
\begin{equation}
r=:\lim_{t\rightarrow \infty }\dfrac{\mathrm{Var}\int_{0}^{t}f^{\prime
}(x_{\tau })d\tau }{t}\geq 0.
\end{equation}

The rate $r=0$ if and only if $\mathrm{Var}(f^{\prime }(x))=0$. Indeed, if $%
\mathrm{Var}(f^{\prime }(x))>0$ then 
\begin{equation}
r(t)=\frac{1}{t}\mathrm{Var(}f^{\prime }(x))\int_{0}^{t}\int_{0}^{t}\mathrm{%
Corr}(f^{\prime }(x_{\tau }),f^{\prime }(x_{s}))d{\tau }ds.  \label{lever:rt}
\end{equation}%
According to Proposition \ref{lever:propo1} in Appendix \ref%
{lever:AuxiliaryStatisticalResult},\textrm{\ }$\mathrm{Corr}(f^{\prime
}(x_{\tau }),f^{\prime }(x_{s}))>0$, so it follows from (\ref{lever:rt})
that 
\begin{equation}
r(t)\geqslant \frac{1}{t}\mathrm{Var}(f^{\prime }(x))\int_{0}^{t}1d{\tau }=%
\mathrm{Var}(f^{\prime }(x))>0.
\end{equation}

Finally, since\textrm{\ }$\mathrm{Var}(f^{\prime }(x))=0$ if and only if $%
f^{\prime }(x)$ is almost surely constant, the investment policy must be
linear if the variance of $u_{t}$ is not to increase with time. QED.

Theorem \textbf{\ref{lever:proposition1} }shows that all linear strategies
eliminate long-run risk. As a natural consequence, the next theorem
predicates optimality of linear strategies with respect to the asymptotic
investment criterion. The idea is to match any non-linear strategy with an
admissible linear strategy that has higher expected return and lower growth
in variance. The matching is possible exactly because all linear strategies
have zero asymptotic growth in variance.

\begin{theorem}
\label{lever:proposition2} For the investor with asymptotic preferences, any 
$\mathbf{D}-$policy is dominated by some linear $\mathbf{D}-$policy.
\end{theorem}

\noindent \textbf{Proof:} Let the non-linear policy be $f(x)$. According to (%
\ref{representation}), in the long run%
\begin{equation}
E(u_{t})=u_{0}-\frac{\sigma ^{2}}{2}E(f^{\prime }(x))t.
\end{equation}%
Take the linear policy $f_{L}(x)=(E(f^{\prime }(x))-\varepsilon )x$ with $%
\varepsilon >0$. For a certain $\varepsilon $ it is admissible. This is
because $|E(f^{\prime }(x))-\varepsilon |<K$ follows from $|f^{\prime
}(x)|\leq K$ everywhere and $|f^{\prime }(x)|<K$ on a set of positive measure%
$,$ which are both true because $f$ is a non-linear $\mathbf{D}-$policy.\
The expectation of the logarithm of wealth under $f_{L}$ is 
\begin{equation}
E(u_{t})=u_{0}-\frac{\sigma ^{2}}{2}E(f^{\prime }(x))t+\frac{\sigma ^{2}}{2}%
\varepsilon t.
\end{equation}%
It is clearly higher than the corresponding expectation for the non-linear
policy. From Theorem \ref{lever:proposition1} we know that the variance of $%
u_{t}$ is asymptotically constant for linear policies and is asymptotically
equivalent to $rt$, where $r>0,$ for non-linear policies. It follows that
for sufficiently large $T$ the linear policy $f_{L}$ will have lower $%
\mathrm{Var}(u_{T})$ than the non-linear policy $f$. Thus $f_{L}$
asymptotically dominates $f$. QED.

It remains to find the optimal policy in the class of linear policies. It
turns out that it is the policy that has the maximal sensitivity to
mispricing. This is intuitively clear because every linear strategy
eliminates the long-run risk, and the policy with the largest sensitivity to
mispricing has the highest expected return. Formally, the following theorem
holds.

\begin{theorem}
\label{lever:proposition3} The optimal linear $\mathbf{D}-$policy is $%
f(x)=-Kx.$
\end{theorem}

\noindent \textbf{Proof: }it is easy to compute 
\begin{eqnarray*}
\lim_{t\rightarrow \infty }\frac{E(u_{t})}{t} &=&\frac{\sigma ^{2}}{2}k, \\
\lim_{t\rightarrow \infty }\frac{\mathrm{Var}(u_{t})}{t} &=&0, \\
U &=&\frac{\sigma ^{2}}{2}k.
\end{eqnarray*}%
Thus, the utility is maximized by the maximal possible $k,$ from which the
theorem follows. QED.

The theorem implies that the arbitrageur should increase the leverage at the
maximal possible rate. In particular, the optimal strategy does not depend
on the risk aversion parameter or properties of the mispricing process. The
intuitive meaning of this conclusion is that the appropriate use of leverage
allows the arbitrageur to eliminate all the long run risk. As the next
section shows, this conclusion will be reversed if the arbitrageur is more
constrained in the use of leverage.

\section{Optimal Threshold Policy}

\label{lever:OptimalThresholdStrategy}

When an arbitrageur uses \textbf{threshold policies} he keeps his finger on
a button that triggers investment while looking at the computer monitor and
waiting for a mispricing. If he observes a mispricing that exceeds a
threshold, $S,$ he pushes the button and a fixed amount, $L,$ is directed to
this opportunity. When the mispricing falls below another threshold, $s,$ he
pushes another button and the position closes. Leverage $L$ never changes
when the position is opened$.$ \textbf{Simple threshold policies} have equal
thresholds: $S=s.$

As was said in the Introduction, arbitrageurs use threshold policies because
they often cannot secure additional funds for positions they already
opened.\ They also favor threshold policies because these policies allow
economizing on transaction costs.

General threshold policies are complicated to analyze. Fortunately, the
following theorem shows that it is sufficient to study simple threshold
policies.

\begin{theorem}
\label{lever:proposition4} Any threshold policy is dominated by a simple
threshold policy.
\end{theorem}

\noindent This theorem is given without proof. Intuitively, for the Markov
process of mispricing the optimal investment policy should not depend on the
history of investing, and the only threshold policies that pass this
selection test are simple threshold policies. Indeed, if $s<S,$ and the
mispricing is between $s$ and $S,$ then the position is open if the
mispricing has fallen from above $S$ but not yet gone below $s,$ and it is
closed if the mispricing has risen from below $s$ but not yet gone above $S.$
It follows that investment under the threshold rule with $s\neq S$ depends
on history of investment and therefore cannot be optimal.

The relevant properties of the simple threshold policies are described in
the next theorem, which uses the following notation:%
\begin{align}
\phi (S)& =\frac{1}{\sqrt{2\pi \Sigma }}\exp \left( -\frac{S^{2}}{2\Sigma }%
\right) , \\
\psi (S)& =\frac{1}{\sqrt{2\pi \Sigma }}\frac{2}{\alpha }\int_{0}^{1}\frac{1%
}{\xi }\left[ \frac{1}{\sqrt{1-\xi ^{2}}}\exp \left( \frac{S^{2}}{\Sigma }%
\frac{\xi }{1+\xi }\right) -1\right] d\xi .
\end{align}%
For $S=0$, the value of $\psi (S)$ can be computed explicitly: $\psi
(0)=2\ln 2/(\alpha \sqrt{2\pi \Sigma }).$

\begin{theorem}
\label{lever:proposition5} For the simple threshold policy with threshold $S$
and leverage~$L$%
\begin{align*}
\lim_{t\rightarrow \infty }\frac{E(u_{t}-u_{0})}{t}& =c_{1}(L,S)\equiv
\sigma ^{2}L\phi (S) \\
\lim_{t\rightarrow \infty }\frac{\mathrm{Var}(u_{t}-u_{0})}{t}&
=c_{2}(L,S)\equiv (\sigma ^{2}L\phi (S))^{2}\psi (S).
\end{align*}%
The investor's utility is $U(L,S)=c_{1}(L,S)-\gamma c_{2}(L,S).$
\end{theorem}

\noindent The proof is relegated to Appendix \ref%
{lever:ProofTheoremThreshold}.

The first step in obtaining the optimal policy from this theorem is to
calculate reduced utility function that depends only on threshold $S:$

\begin{corollary}
For a fixed threshold $S$ the optimal leverage is 
\begin{equation*}
L(S)=\frac{1}{4\gamma \alpha \Sigma \phi (S)\psi (S)}
\end{equation*}%
and the corresponding utility is 
\begin{equation*}
U(S)=\frac{1}{4\gamma \psi (S)}.
\end{equation*}
\end{corollary}

The functions $L(S)$ and $U(S)$ are illustrated in Figures~\ref%
{lever:figure1a} and \ref{lever:figure1b}. We can see that higher long-run
variance $\Sigma $ leads to an increase in both leverage $L$ and utility $U$%
. Higher persistence \ of the process does not change optimal leverage but
decreases utility.

\begin{figure}[tbp]
\framebox[\textwidth]{[Figure \ref{lever:figure1a} here]}
\caption{Plots of Optimal Leverage $L$ as Function of Threshold $S$}
\label{lever:figure1a}
\end{figure}

\begin{figure}[tbp]
\framebox[\textwidth]{[Figure \ref{lever:figure1b} here]}
\caption{Plots of Optimal Utility $U$ as Function of Threshold $S$}
\label{lever:figure1b}
\end{figure}

The function $\psi (S)$ is increasing in $S^{2},$ and consequently the
maximal utility is reached at $S=0.$ The optimal threshold policy is given
in the following theorem.

\begin{theorem}
Utility is maximized for $S=0$ and $L=\pi /(4\gamma \ln 2)$. The optimal
utility is $U=\alpha \sqrt{2\pi \Sigma }/(8\gamma \ln 2).$
\end{theorem}

Predictably, the utility is higher when the convergence is fast ($\alpha $
is high) and the arbitrage opportunity is large ($\Sigma $ is high). Not so
predictable is that the optimal leverage does not depend on the parameters
of the process: This leverage optimally balances risk and return for every
Ornstein-Uhlenbeck process. What is most important, however, is that the
optimal leverage depends on the parameter of the risk aversion $\gamma $.
The higher $\gamma $ is, the lower the amount is that the arbitrageur is
willing to commit to the arbitrage opportunity: The arbitrageur that uses
only threshold strategies is unable to remove the long-run risk and must
adjust his behavior.

\section{Empirical Application}

\label{lever:EmpiricalApplication} This section studies convergent trading
in the context of WEBS, which are shares of open-end mutual funds that
replicate the price performance of foreign stock market indexes. WEBS trade
on a stock exchange like ordinary stock, and their managers try to keep the
fund price close to the net asset value (NAV) of their underlying stocks.
They buy back shares if the price is less then NAV and issue additional
shares if the price is greater than NAV.

As in previous sections, I assume that the investor can hedge the risk of
the underlying portfolio. Indeed, trading index futures provides a very good
hedge of country exposure. Absense of the perfect hedge limits the
implications of my analysis.

I use price and NAV daily data for WEBS that track indices of Australia,
Austria, Belgium, Canada, France, Germany, Hong Kong, Italy, Japan,
Malaysia, Mexico, Netherlands, Singapore, Spain, Sweden, Switzerland, and
the United Kingdom. This is a total of 17 countries. The data start in March
of 1996 and end in August of 2000, which gives around 1000 datapoints for
each country.

The mispricing factor $x_t$ is computed as the logarithm of the ratio of the
price to NAV. Some summary statistics for $x_t$ are given in Table~\ref%
{lever:table1}. I model the dynamic of the mispricing factor as an AR(1)
process, which is the discrete-time counterpart to the Ornstein-Uhlenbeck
process: 
\begin{equation}
x_t=\beta x_{t-1} +\sigma \epsilon_{t},
\end{equation}

\begin{table}[tbp]
\framebox[\textwidth]{[Table \ref{lever:table1} here]}
\caption{Summary Statistics of Mispricing Factor}
\label{lever:table1}
\end{table}

The results of estimation of this process are summarized in Table~\ref%
{lever:table2}. They show that $\beta$ is around $0.5$, and $\sigma$ is
around $0.01$. Durbin-Watson statistic shows that the AR(1) process is a
reasonably good approximation to the true mispricng process.

\begin{table}[tbp]
\framebox[\textwidth]{[Table \ref{lever:table2} here]}
\caption{Results of Estimation of Mispricing Factor Process}
\label{lever:table2}
\end{table}

Our first goal is to get an estimate of the order of the coefficient $k$ in
the optimal linear strategy. Let us use the following estimates of the order
of parameters: $\sigma^2\sim 10^{-4}$, $\alpha\sim 0.5$. From Theorem \ref%
{lever:proposition3}, the long-run variance of the logarithm of wealth is $%
5\cdot10^{-9}k^2,$ and the expected change in the logarithm of wealth is $%
5\cdot 10^{-5}k$ per day. The average daily change in the logarithm of the S$%
\&$P500 index has been $8\cdot 10^{-4}\sim 10^{-3}$ in the last five years.
To get this return by convergent trading, the sensitivity $k$ of the linear
strategy would have to be set at $20.$ The corresponding asymptotic variance
of the logarithm of wealth would then be $2\cdot 10^{-8}$. This is
considerably smaller than the variance of the deviation of the logarithm of
the $S\&P500$ index from its linear trend, which for the 5-year period from
August 1995 to August 2000 can be estimated at about $3\cdot 10^{-3}$.

These computations suggest that this market could not be efficient if
transactions costs were absent. The goal of the following Monte-Carlo
simulations is to analyze properties of threshold strategies and to find the
optimal threshold strategy in situations with transaction costs. I assume
that the transaction cost $c$ is $0.25\%$, and that the true process of
discounts is AR(1) with $\beta=0.3$ and $\sigma=10^{-2}$.

The simulations were organized as follows. One hundred realizations of the
mispricing factor process were generated. Each realization had $1250$
datapoints. For each realization I simulated the process of investing
according to a strategy from a finites set of threshold strategies. Thus,
each pair of a strategy and a realization of the mispricing factor process
resulted in a realization of wealth process.

For each realization of wealth process, I calculated the average daily
return as a difference between logarithm of final wealth and logarithm of
initial wealth divided by the length of the realization. After that, I
calculated the mean and the variance of this statistic over all realizations
of the wealth process that corresponded to a particular strategy. Thus, as a
final product I had a function that mapped each strategy into the
expectation and variance of the average daily return.

I used the following set of strategies. A rise in the mispricing over a
threshold $S$ triggers the opening of the position. When the mispicing
returns to the region below $s\in \lbrack 0,S]$, the arbitrageur closes his
position and waits for a new trigger signal. The threshold $S$ was chosen in
the range from $0.5\%$ $\ $to $2\%$. The threshold $s$ was chosen in the
range from $S+c$ to $0\%$.

\begin{figure}[tbp]
\framebox[\textwidth]{[Figure \ref{lever:figure2} here]}
\caption{Contour Graph of Mean of Average Daily Returns}
\label{lever:figure2}
\end{figure}

\begin{figure}[tbp]
\framebox[\textwidth]{[Figure \ref{lever:figure3} here]}
\caption{Contour Graph of Standard Deviation of Average Daily Returns}
\label{lever:figure3}
\end{figure}

\begin{figure}[tbp]
\framebox[\textwidth]{[Figure \ref{lever:figure4} here]}
\caption{Contour Graph of Sharpe Ratio}
\label{lever:figure4}
\end{figure}

Figure~\ref{lever:figure2} is a contour graph of the mean of the average
daily return to a strategy. On the vertical axis of the graph is the low
threshold $S$ and on the horizontal axis is the difference between high and
low thresholds $S-s$. They are denominated in percentage terms. The lines on
the contour graph correspond to the strategies that have the same mean of
average daily return.

This graph suggests that the average return is maximized for $s=0$ and $S$
set to some $S_0>c$. Therefore, if the investor wants to maximize average
return he should invest only if the mispricing exceeds the transaction cost
by some markup.

Figure~\ref{lever:figure3} is a contour graph of the standard deviation of
the average daily returns. We can conclude from it that the variance is
increasing with an increase in the absolute value of the low threshold and
with an increase in the difference between the thresholds.

Figure~\ref{lever:figure4} is a countour graph of the ratio of the mean of
average daily return to its standard deviation. It suggests that the ratio
is maximized for the strategy that sets $s=0$ and $S=c$. Thus an investor
who uses only threshold strategies and who wants to maximize the Sharpe
ratio should invest immediately when the mispricing exceeds his transaction
costs.

\section{Discussion}

Section \ref{lever:OptimalDifferentiableStrategy} shows that in the class of
differentiable policies with bounded derivative the optimal policy is linear
in the mispricing and the coefficient in the linear relationship is the
highest possible. The optimal strategy in this case does not depend on the
risk-aversion of the arbitrageur, and all the long-run risk can be
eliminated.

In contrast, according to the results of Section \ref%
{lever:OptimalThresholdStrategy}, if only threshold policies are available
then the long-run risk is unavoidable and the investment is inversely
proportional to risk aversion. This conclusion is consistent with the belief
that arbitrageurs are typically risk-neutral. The suggested reason for this
belief is that the constraints on flexibility of changes in leverage make
the convergence trading risky even in the long run.

\newpage

\appendix

\section{Auxiliary Statistical Result}

\label{lever:AuxiliaryStatisticalResult}

Suppose that $x$ and $y$ are jointly Gaussian random variables,\textrm{\ }$%
\mathrm{Var}(x)=\mathrm{Var}(y)=1,$ and $\mathrm{Cov}(x,y)=\beta $.

\begin{proposition}
\label{lever:propo1} If $f\in \mathbf{D}$ and $\mathrm{Var}(f(x))=1$, then $%
\mathrm{Cov}(f(x),f(y))\in \lbrack 0,\beta ]$.
\end{proposition}

\noindent \textbf{Proof:} Assume without loss of generality that $Ef(x)=0$.
Since \ Hermite polynomials are complete in the class of $\mathbf{D}-$%
policies, we can use them to approximate $f$. Then the assertion of
Proposition \ref{lever:propo1} follows from Proposition \ref{lever:propo2}.

\begin{proposition}
\label{lever:propo2} If $f$ is a polynomial of degree $N$, $Ef(x)=0$ and $%
\mathrm{Var}(f(x))=1$, then $\mathrm{Cov}(f(x),f(y))\in \lbrack \beta
^{N},\beta ].$ The maximum and minimum are achieved by $f(x)=x$ and $%
f(x)=H_{N}(x)$, respectively, where $H_{N}(x)$ is the Hermite polynomial of
degree $N.$
\end{proposition}

\noindent \textbf{Proof:} Represent $f(x)$ as a sum of Hermite polynomials: 
\begin{equation}
f(x)=\sum_{1}^{N}a_{k}H_{k}(x),
\end{equation}%
where by definition 
\begin{equation}
H_{k}(x)=\exp \left( \frac{x^{2}}{2}\right) \frac{(-1)^{k}}{k!}\frac{d^{k}}{%
(dx)^{k}}\left[ \exp \left( -\frac{x^{2}}{2}\right) \right] .
\end{equation}%
Hermite polynomials form an orthonormal system with respect to the Gaussian
kernel and possess the following useful property: 
\begin{equation}
\mathrm{Cov}(H_{i}(x),H_{j}(y))=\beta ^{i}\delta _{ij}.
\end{equation}%
Using this property and orthonormality, we can write 
\begin{align}
\mathrm{Cov}(f(x),f(y))& =\sum_{1}^{N}a_{k}^{2}\beta ^{k}  \label{lever:cov}
\\
\mathrm{Var}(f(x))& =\sum_{1}^{N}a_{k}^{2}.  \label{lever:var}
\end{align}%
From (\ref{lever:cov}) and (\ref{lever:var}), the maximum of $\mathrm{Cov}%
(f(x),f(y))$ is $\beta $ and it is achieved by $f(x)=H_{1}(x)=x$. The
minimum is $\beta ^{N}$ and it is achieved by $f(x)=H_{N}(x)$. QED.

\section{Proof of Theorem \ref{lever:proposition5}}

\noindent \label{lever:ProofTheoremThreshold}\textbf{Proof:} By definition,
the threshold policy is 
\begin{equation}
f(x)=%
\begin{cases}
-\mathrm{sign}(x)L & \text{ if $|x|\geq S$}, \\ 
0 & \text{ if $|x|<S$}.%
\end{cases}%
\end{equation}%
The generalized Ito formula gives 
\begin{equation}
u_{t}=u_{0}+g(x_{t})-g(x_{0})+\frac{\sigma ^{2}}{2}2L\int_{0}^{t}\delta
_{S}(x_{\tau })d\tau ,  \label{lever:LocalTime}
\end{equation}%
where $\delta _{S}$ is the Dirac delta-function and 
\begin{equation}
g(\xi )=:\int_{0}^{\xi }f(\zeta )d\zeta .
\end{equation}%
The intuition behind this representation is simple: The investor increases
his wealth only when he triggers the policy. The number of times the policy
is triggered is stochastic and measured by the integral of the delta
function. The profit earned at each occasion is proportional to the product
of local volatility $\sigma ^{2}$ and leverage $L.$ Finally, there is a
dependence of wealth on initial and final conditions which is captured by $%
g(x_{t})-g(x_{0}).$

Since $g(x_{t})$ does not grow with time, the arbitrageur's utility depends
only on the moments of the integral of the delta function:%
\begin{equation}
\int_{0}^{t}\delta _{S}(x_{\tau })\,d\tau .
\end{equation}

The first step in the computation of the moments is calculating the
expectation and covariance function of the generalized stochastic process $%
\delta _{t,S}=:\delta _{S}(x_{t})$. The joint density of $x_{t_{1}}$ and $%
x_{t_{2}}$ is 
\begin{equation}
p(x_{1},x_{2})=\frac{1}{2\pi \Sigma \sqrt{1-a(\tau )^{2}}}\exp \left\{ -%
\frac{1}{2\Sigma }%
\begin{pmatrix}
x_{1} \\ 
x_{2}%
\end{pmatrix}%
^{\prime }%
\begin{pmatrix}
1 & a(\tau ) \\ 
a(\tau ) & 1%
\end{pmatrix}%
^{-1}%
\begin{pmatrix}
x_{1} \\ 
x_{2}%
\end{pmatrix}%
\right\} ,
\end{equation}%
where $\tau =t_{2}-t_{1}$ and $a(\tau )=e^{-\alpha |\tau |}$. The delta
function can be approximated by $\frac{1}{\Delta }\chi _{\lbrack S,S+\Delta
]}$, where $\chi _{A}$ denotes the characteristic function of set $A$ and $%
\Delta $ limits to 0. Then, computing two first moments for $\chi _{\lbrack
S,S+\Delta ]}(x_{t})$ and taking the limit $\Delta \rightarrow 0$ give the
following formulas: 
\begin{align}
E(\delta _{t_{1},S})& =\frac{1}{\sqrt{2\pi \Sigma }}\exp \left( -\frac{S^{2}%
}{2\Sigma }\right) \text{ and}  \label{lever:DeltaMean} \\
E(\delta _{t_{1},S}\delta _{t_{2},S})& =\frac{1}{2\pi \Sigma \sqrt{1-a(\tau
)^{2}}}\exp \left( -\frac{S^{2}}{\Sigma }\frac{1}{1+a(\tau )}\right) .
\label{lever:DeltaVar}
\end{align}%
For example, equality (\ref{lever:DeltaVar}) can be seen from the following
calculation:

\begin{equation}
\begin{split}
& E(\delta _{t_{1},S}\delta _{t_{2},S})=\lim_{\Delta \rightarrow 0}\frac{1}{%
\Delta ^{2}}\int_{-\infty }^{\infty }\int_{-\infty }^{\infty }\chi _{\lbrack
S,S+\Delta ]}(\xi _{1})\chi _{\lbrack S,S+\Delta ]}(\xi _{2})p(\xi _{1},\xi
_{2})d\xi _{1}d\xi _{2} \\
& =\lim_{\Delta \rightarrow 0}\frac{1}{2\pi \Sigma \sqrt{1-a(\tau )^{2}}}%
\frac{1}{\Delta ^{2}}\int_{S}^{S+\Delta }\int_{S}^{S+\Delta }\exp \left[ -%
\frac{\xi _{1}^{2}+\xi _{2}^{2}-2a(\tau )\xi _{1}\xi _{2}}{2\Sigma (1-a(\tau
)^{2})}\right] d\xi _{1}d\xi _{2} \\
& =\frac{1}{2\pi \Sigma \sqrt{1-a(\tau )^{2}}}\exp \left( -\frac{S^{2}}{%
\Sigma }\frac{1}{1+a(\tau )}\right) .
\end{split}%
\end{equation}

From (\ref{lever:DeltaMean}) and (\ref{lever:DeltaVar}) it follows that 
\begin{equation}
\mathrm{Cov}(\delta _{t_{1},S},\delta _{t_{2},S})=\frac{1}{2\pi \Sigma \sqrt{%
1-a(\tau )^{2}}}\exp \left( -\frac{S^{2}}{\Sigma }\frac{1}{1+a(\tau )}%
\right) -\frac{1}{2\pi \Sigma }\exp \left( -\frac{S^{2}}{\Sigma }\right) .
\end{equation}%
Since $\mathrm{Cov}(\delta _{t_{1},S},\delta _{t_{2},S})$ depends only on $%
\tau =t_{2}-t_{1},$ it can be denoted by $\vartheta (\tau ,S).$ Next,%
\begin{equation}
\begin{split}
\mathrm{Var}\left( \int_{0}^{T}\delta _{t,S}\,dt\right) &
=2\int_{0}^{T}dt_{1}\int_{t_{1}}^{T}\mathrm{Cov}(\delta _{t_{1},S},\delta
_{t_{2},S})dt_{2} \\
\text{(substituting $\tau =t_{2}-t_{1}$)}&
=2\int_{0}^{T}dt_{1}\int_{0}^{T-t_{1}}\vartheta (\tau ,S)d\tau \\
\text{(changing order of integration)}& =2\int_{0}^{T}(T-\tau )\vartheta
(\tau ,S)d\tau .
\end{split}%
\end{equation}

Since $\vartheta (\tau ,S)=O(e^{-c\tau })$ for a positive $c$ and $\tau
\rightarrow \infty $, and $\vartheta (\tau ,S)$ is integrable around $\tau
=0,$ it follows%
\begin{align*}
& \int_{0}^{T}\vartheta (\tau ,S)d\tau \rightarrow \mathrm{const}\text{ and}
\\
& \int_{0}^{T}\tau \vartheta (\tau ,S)d\tau \rightarrow \mathrm{const}
\end{align*}%
as $T\rightarrow \infty $.

So, for large~$T$ 
\begin{equation}
\begin{split}
& \mathrm{Var}\left( \int_{0}^{T}\delta _{t,S}\,dt\right) \sim
2T\int_{0}^{\infty }\vartheta (\tau ,S)d\tau \\
& =2T\frac{1}{2\pi \Sigma }\exp \left( -\frac{S^{2}}{\Sigma }\right)
\int_{0}^{\infty }\left( \frac{1}{\sqrt{1-a(\tau )^{2}}}\exp \left( \frac{%
S^{2}}{\Sigma }\frac{a(\tau )}{1+a(\tau )}\right) -1\right) \,d\tau \\
& \text{(substituting $\tau =-\alpha ^{-1}\ln \xi $)} \\
& =2T\frac{1}{2\pi \Sigma }\exp \left( -\frac{S^{2}}{\Sigma }\right) \frac{1%
}{\alpha }\int_{0}^{1}\frac{1}{\xi }\left( \frac{1}{\sqrt{1-\xi ^{2}}}\exp
\left( \frac{S^{2}}{\Sigma }\frac{\xi }{1+\xi }\right) -1\right) \,d\xi \\
& =\phi (S)^{2}\psi (S)T.
\end{split}%
\end{equation}%
This implies all the assertions of the theorem. QED.

\newpage

\bibliographystyle{CHICAGO}
\bibliography{comtest}

\end{document}